# On the Computation of Hilbert Bases and Extreme Rays of Cones

Raymond Hemmecke,
Department of Mathematics, University of California, Davis,
One Shields Avenue, Davis, CA 95616, USA
raymond@hemmecke.de

October 24, 2018

## Abstract

In this paper we present a novel project-and-lift approach to compute the set of minimal generators of the semigroup $(\Lambda \cap \mathbb{R}^n_+, +)$ for lattices $\Lambda \subseteq \mathbb{Z}^n$. This problem class includes the computation of Hilbert bases of cones $\{z : Az = 0, z \in \mathbb{R}^n_+\}$ for integer matrices $A$. A similar approach can be used to compute only the extreme rays of such cones. Finally, some combinatorial applications and computational experience are presented.

## 1 Introduction

Lattice points in polyhedral cones arise as interesting objects in many branches of mathematics as, for example, combinatorics, integer programming, computational algebra, or topology. Often one is interested in a finite subset of lattice points in the cone that generate all other lattice points in the cone as non-negative integer linear combinations.

**Definition 1.1** *Let $\mathcal{C} \subseteq \mathbb{R}^n$ be a polyhedral cone with rational generators and let $\Lambda \subseteq \mathbb{Z}^n$ be a lattice. Then we call a finite set $H = \{h_1, \ldots, h_t\} \subseteq \Lambda \cap \mathcal{C}$ a generating set of the monoid $(\Lambda \cap \mathcal{C}, +)$ if for every $z \in \Lambda \cap \mathcal{C}$ there are non-negative integral multipliers $\lambda_1, \ldots, \lambda_t$ such that $z = \sum_{i=1}^t \lambda_i h_i$.*

Note that for every pointed rational polyhedral cone $\mathcal{C} \subseteq \mathbb{R}^n$ and every lattice $\Lambda \subseteq \mathbb{Z}^n$ there exists a unique generating set of the monoid $(\Lambda \cap \mathcal{C}, +)$ that is minimal with respect to inclusion ([12], Chapter 7). For $\Lambda = \mathbb{Z}^n$ this generating set is also called a Hilbert basis [13]. It is in general a hard problem to compute Hilbert bases. In this paper we will give a novel project-and-lift algorithm to compute minimal generating sets for $\mathcal{C} = \mathbb{R}^n_+$ and arbitrary lattices $\Lambda \subseteq \mathbb{Z}^n$. One prominent special case, for which many algorithms have been proposed (e. g. [2, 9, 10, 11, 12, 14]), is the case of $\Lambda = \ker_\mathbb{Z}(A)$ for matrices $A \in \mathbb{Z}^{d \times n}$, where $\ker_\mathbb{Z}(A) = \{z : Az = 0, z \in \mathbb{Z}^n\}$ denotes the integral kernel of $A$. In these algorithms, the Hilbert basis is extracted from a (usually much) bigger superset, in some cases even additional variables have to be introduced. The advantage of the algorithm that we propose here is that we do not need additional variables, that we may throw away unnecessary vectors in early and intermediate stages of the algorithm, and that we therefore arrive at a far smaller superset from which the set of minimal generators of the monoid can be extracted. A similar project-and-lift approach leads to a novel algorithm to compute the extremal rays of cones given as $\{z : Az = 0, z \in \mathbb{R}^n_+\}$, see Section 4. Note that one may apply this approach also for the computation of truncated Graver bases which improves the algorithm given in [8]. Some implementations of algorithms to compute extreme rays of cones (e. g. cdd [7], porta [4]) and to compute Hilbert bases (e. g. NORMALIZ [2]) are freely available over the internet. The latter code, however, needs a generating set of the cone $\{z : Az = 0, z \in \mathbb{R}^n_+\}$ as input. Thus, no direct computational comparison between NORMALIZ and the algorithm presented below is possible. The main objective of this paper is an algorithmic solution to the following problem.

**Problem 1.2** *Given $\{p_1, \ldots, p_s\} \subseteq \mathbb{Z}^n$ that generate a lattice $\Lambda \subseteq \mathbb{Z}^n$ over $\mathbb{Z}$, compute the unique minimal generating set $H(\Lambda)$ of the monoid $(\Lambda \cap \mathbb{R}^n_+, +)$ and the set $R(\Lambda)$ of extremal rays of the pointed rational cone $\mathcal{C}_\Lambda = \{z : z = \sum_{i=1}^s \lambda_i p_i, \lambda_1, \ldots, \lambda_s \in \mathbb{R}, z \in \mathbb{R}^n_+\}$.*

In Section 3 we present a selection strategy for the S-vector that has to be considered next in the Completion Algorithm 2.3 to compute the minimal generators of the monoid. This strategy leads to a tremendous speed-up of our project-and-lift algorithm. In Section 5 we collect some applications of Hilbert bases and extreme ray computations. These include, for example, the computation of dual cones. The extreme rays of the dual cone of a cone $\mathcal{C}$, however, are exactly the normal vectors of the facets of $\mathcal{C}$. Finally, in Section 6, we report on an implementation, MLP, and some computational experience while working on the problem of counting magic arrays [1].

## 2 Computation of Minimal Generators

In order to state the algorithm below we need to introduce some useful notation.

**Definition 2.1** *For any $m \geq j$ let $\pi_j^m : \mathbb{R}^m \to \mathbb{R}^j$ be the projection onto the first $j$ coordinates. Moreover, let $K_j := \{\pi_j^n(v) : v \in \Lambda\}$, $K_j^+ := K_j \cap (\mathbb{R}^{j-1}_+ \times \mathbb{R}_+)$, $K_j^- := K_j \cap (\mathbb{R}^{j-1}_+ \times \mathbb{R}_-)$.*

Note that $K_j$ is a sublattice of $\mathbb{Z}^j$ and that $K_j^+$ and $K_j^-$ are intersections of a lattice, $K_j$, with pointed rational cones, that is, $(K_j^+, +)$ and $(K_j^-, +)$ are monoids. Let $H_j^+$ and $H_j^-$ denote their unique inclusion minimal generating sets. Clearly, $H(\Lambda) = H_n^+$ since $K_n^+ = \Lambda \cap \mathbb{R}^n_+$.
We may assume that the generators $\{p_1, \ldots, p_s\} \subseteq \mathbb{Z}^n$ of $\Lambda$ have the following structure, which is clearly achievable by suitable elementary integral row operations and permutations of columns:

$$\begin{array}{rcrrcccccc}
p_1 & = & ( & p_{1,1}, & p_{1,2}, & \ldots, & \ldots, & p_{1,s}, & \ldots, & p_{1,n} \ ), \\
p_2 & = & ( & 0, & p_{2,2}, & \ldots, & \ldots, & p_{2,s}, & \ldots, & p_{2,n} \ ), \\
p_3 & = & ( & 0, & 0, & p_{3,3} & \ldots, & p_{3,s}, & \ldots, & p_{3,n} \ ), \\
& \vdots & ( & \vdots, & \vdots, & \ldots, & \ddots, & \vdots, & \ldots, & \vdots \ ), \\
p_s & = & ( & 0, & 0, & \ldots, & 0, & p_{s,s}, & \ldots, & p_{s,n} \ ),
\end{array}$$

with $p_{i,i} > 0$, $i = 1, \ldots, s$.
In what follows, we will start with $H_1^+ = \{(p_{1,1})\}$, which is easily computed, and compute $H_2^+$, $H_3^+$, …, $H_n^+$ inductively. For each step $H_j^+ \to H_{j+1}^+$ we will employ a completion procedure [3]. The input sets for $j \geq s$ and $j < s$ differ slightly due to the fact that for $j < s$ the map $\pi_j^n$ is not injective. In particular, there are non-zero vectors $v \in \Lambda$ with $\pi_j^n(v) = 0$ if $j < s$.
As $H_j^+ \subseteq K_j$, we know that for every $h \in H_j^+$ there is at least one vector $v \in \Lambda$ with $h = \pi_j^n(v)$. Thus, for every $h \in H_j^+$ there is some $h' \in \mathbb{Z}$ such that $(h, h') \in K_{j+1}$. If there are several choices for $h'$, which can only happen if $j < s$, we will choose the smallest non-negative integer $h'$ such that $(h, h') \in K_{j+1}$ to form the input set $F$ below. For $j < s$ we have to add the vectors $\pi_{j+1}^n(p_{j+1})$ and $-\pi_{j+1}^n(p_{j+1})$ to the input set $F$.

**Lemma 2.2** *The set*

$$F = \begin{cases} \bigcup_{h \in H_j^+} \{(h, h') : (h, h') \in K_{j+1}\} & \text{if } j \geq s, \text{ and} \\ \bigcup_{h \in H_j^+} \{(h, h') : (h, h') \in K_{j+1}\} \cup \{\pi_{j+1}^n(p_{j+1}), -\pi_{j+1}^n(p_{j+1})\} & \text{if } j < s, \end{cases}$$

*generates $K_{j+1}^+ \cup K_{j+1}^-$ over $\mathbb{Z}_+$, that is, every $z \in K_{j+1}^+ \cup K_{j+1}^-$ can be written as a non-negative integer linear combination of elements of $F$.*

**Proof.** Let $\pi_{j+1}^n(v) \in K_{j+1}^+ \cup K_{j+1}^-$, where $v \in \Lambda$. Since $H_j^+$ is a generating set for the monoid $(K_j^+, +)$, we can write $\pi_j^n(v) = \sum \alpha_i \pi_j^n(v_i)$ for some $\pi_j^n(v_i) \in H_j^+$ and $\alpha_i > 0$. Therefore, $\pi_{j+1}^n(v) - \sum \alpha_i \pi_{j+1}^n(v_i)$ has zeros in the first $j$ components and must be the zero vector if $j \geq s$, or an integral multiple of $\pi_{j+1}^n(p_{j+1})$ if $j < s$. This concludes the proof. □

In what follows, dashed variables, like $h'$, at the end of a vector, like $(h, h')$, always represent integer numbers.

**Algorithm 2.3** *(Algorithm to Compute $H_{j+1}^+ \cup H_{j+1}^-$)*

<u>Input:</u> $F = \begin{cases} \bigcup_{h \in H_j^+} \{(h, h') : (h, h') \in K_{j+1}\} & \text{if } j \geq s, \text{ and} \\ \bigcup_{h \in H_j^+} \{(h, h') : (h, h') \in K_{j+1}\} \cup \{\pi_{j+1}^n(p_{j+1}), -\pi_{j+1}^n(p_{j+1})\} & \text{if } j < s \end{cases}$

<u>Output:</u> a set $G$ which contains $H_{j+1}^+ \cup H_{j+1}^-$

$G := F$
$C := \bigcup_{f, g \in G} \text{S-vectors}(f, g)$
<u>while</u> $C \neq \emptyset$ <u>do</u>
    $s :=$ *an element in* $C$
    $C := C \setminus \{s\}$
    $f := \text{normalForm}(s, G)$
    <u>if</u> $f \neq 0$ <u>then</u>
        $C := C \cup \bigcup_{g \in G} \text{S-vectors}(f, g)$
        $G := G \cup \{f\}$
<u>return</u> $G$.

It remains to define the set S-vectors$(f, g)$ and the function normalForm$(s, G)$.

$$\text{S-vectors}((v, v'), (w, w')) := \begin{cases} \{(v + w, v' + w')\} & \text{if } v'w' < 0, \\ \emptyset & \text{otherwise.} \end{cases}$$

In case that one is interested only in those elements of $H(\Lambda)$ whose components lie below certain bounds $u \in (\mathbb{Z}_+ \cup \{\infty\})^n$, that is $H(\Lambda) \cap \{z : z \leq u, z \in \mathbb{Z}_+^n\}$, these upper bounds $u$ can easily be used to speed-up the computation by considering only the following set of S-vectors:

$$\text{S-vectors}((v, v'), (w, w')) := \begin{cases} \{(v + w, v' + w')\} & \text{if } v'w' < 0 \text{ and } v + w \leq u, \\ \emptyset & \text{otherwise.} \end{cases}$$

This choice will be justified below after the proof of Lemma 2.5.

Behind the function normalForm$(s, G)$ there is the following algorithm, wherein $u \sqsubseteq_{j+1} v$ if and only if $u^{(i)} \leq v^{(i)}$, $i = 1, \ldots, j$, $|u^{(j+1)}| \leq |v^{(j+1)}|$, and $u^{(j+1)}v^{(j+1)} \geq 0$. Note that in case of $u \sqsubseteq_{j+1} v$, $u$ and $v$ lie in the same orthant of $\mathbb{R}^{j+1}$.

**Algorithm 2.4** *(Normal Form Algorithm)*
<u>Input:</u> *a vector $s$, a set $G$ of vectors*
<u>Output:</u> *a normal form of $s$ with respect to $G$*

<u>while</u> *there is some $g \in G$ such that $g \sqsubseteq_{j+1} s$ <u>do</u>*
$\qquad \alpha := \min\{\lfloor s^{(i)}/g^{(i)} \rfloor : i = 1, \ldots, j+1, g^{(i)} \neq 0\}$
$\qquad s := s - \alpha g$
<u>return</u> *s*

Note that the first $j$ components of $s$ always remain non-negative during the normal form algorithm. Moreover, the $(j+1)$st component of $s$ either remains non-negative or non-positive. Since $\|s - \alpha g\|_1 < \|s\|_1$, the normal form algorithm always terminates.

**Lemma 2.5** *Algorithm 2.3 terminates and returns a set containing $H_{j+1}^+ \cup H_{j+1}^-$.*

**Proof.** First, let us show termination of the algorithm. For this define for any $r \in \mathbb{R}$ the symbols $r^+ := \max\{r, 0\}$ and $r^- := \max\{-r, 0\}$. Next, consider the sequence $G \setminus F = \{(g_1, g_1'), \ldots\}$ as it is generated by Algorithm 2.3. This sequence fulfills $(g_i, g_i') \not\sqsubseteq_{j+1} (g_k, g_k')$ whenever $i < k$. Therefore, $(g_i, (g_i')^+, (g_i')^-) \not\leq (g_k, (g_k')^+, (g_k')^-)$ whenever $i < k$. Applying the Gordan-Dickson Lemma (see for example [5]) to the sequence $\{(g_1, (g_1')^+, (g_1')^-), \ldots\} \subseteq \mathbb{Z}_+^{j+2}$, we conclude that this sequence must be finite and thus Algorithm 2.3 terminates. It remains to show correctness of the algorithm. By $G$ denote the set that is returned by Algorithm 2.3 and let $(z, z') \in H_{j+1}^+ \cup H_{j+1}^-$. By Lemma 2.2 and $F \subseteq G$, we can write $(z, z')$ as a finite positive integer linear combination $(z, z') = \sum \alpha_i (v_i, v_i')$ for some $\alpha_i \in \mathbb{Z}_{>0}$ and vectors $(v_i, v_i') \in G$ with $0 \leq v_i \leq z$ for all $i$. From the set of all such positive integer linear combinations $\sum \alpha_i (v_i, v_i')$ choose one such that $\sum \alpha_i |v_i'|$ is minimal. Note that, by the triangle inequality, $\sum \alpha_i |v_i'| \geq |z'|$ with equality if and only if $v_i' \sqsubseteq_1 z'$ for all $i$, that is, if and only if $z'$ and all $v_i'$ are either both non-negative or both non-positive.

If $\sum \alpha_i |v_i'| = |z'|$, then we get $(v_i, v_i') \sqsubseteq_{j+1} (z, z')$ for all $i$. Since $(z, z') \in H_{j+1}^+ \cup H_{j+1}^-$, the representation $(z, z') = \sum \alpha_i (v_i, v_i')$, $\alpha_i \in \mathbb{Z}_{>0}$ must be trivial, that is, $(z, z') = 1 \cdot (z, z') \in G$, and consequently, $(z, z') \in G$ and nothing is left to prove.

Hence, we will assume on the contrary that $\sum \alpha_i |v_i'| > |z'|$ holds. Therefore, there must exist $(v_{i_1}, v_{i_1}'), (v_{i_2}, v_{i_2}')$ such that $v_{i_1}' v_{i_2}' < 0$. The sum $(v_{i_1}, v_{i_1}') + (v_{i_2}, v_{i_2}')$ was reduced to 0 during the run of Algorithm 2.3 which gives a integer linear combination $(v_{i_1}, v_{i_1}') + (v_{i_2}, v_{i_2}') = \sum \beta_k (w_k, w_k')$ for some positive integers $\beta_k$ and some $(w_k, w_k') \in G$.

Moreover, $\beta_k (w_k, w_k') \sqsubseteq_{j+1} (v_{i_1}, v_{i_1}') + (v_{i_2}, v_{i_2}')$ for all $k$, implying that $0 \leq w_k \leq z$ and that

$$\sum \beta_k |w_k'| = |\sum \beta_k w_k'| = |v_{i_1}' + v_{i_2}'| < |v_{i_1}'| + |v_{i_2}'|$$

holds. But then $(z, z') = \sum \beta_k (w_k, w_k') + (\alpha_{i_1} - 1)(v_{i_1}, v_{i_1}') + (\alpha_{i_2} - 1)(v_{i_2}, v_{i_2}') + \sum_{i \neq i_1, i_2} \alpha_i (v_i, v_i')$ contradicts the minimality of $\sum \alpha_i |v_i'|$. We conclude $\sum \alpha_i |v_i'| = |z'|$ and the claim follows. $\square$

Note that if there are upper bounds $u$ on $z$, the relations $v_{i_1} + v_{i_2} \leq z \leq u$ always holds in the above proof. This justifies the above mentioned restriction of necessary S-vectors that need to be considered if upper bounds $u$ on the variables are given.

# 3  Improving the Project-and-Lift Algorithm

In this section we improve Algorithm 2.3 such that we need not reduce the S-vectors from $C$ with respect to $G$. Instead, it suffices to check each S-vector only for reducibility with respect to $G$. Moreover, the set $C$ of S-vectors that still have to be considered for reducibility need not be stored in memory explicitly. Both properties speed-up the computation, save computer memory and thus, allow us to solve bigger problems as by application of Algorithm 2.3 alone.

The major idea of the improved algorithm is to consider the elements $(v, v')$ in $C$ by increasing norm $\|v\|_1 = \sum_{i=1}^{j} v^{(i)}$. To this end, let $\mathcal{G}_i := \{(v, v') \in H_{j+1}^+ \cup H_{j+1}^- : \|v\|_1 = i\}$. Suppose we have already computed $\mathcal{G}_0, \mathcal{G}_1, \ldots, \mathcal{G}_k$, for some $k \geq 1$. Then the following observation allows us to avoid all the reduction steps in order to compute $\mathcal{G}_{k+1}$.

**Lemma 3.1** *Every vector $(z, z') \in K_{j+1}^+ \cup K_{j+1}^-$ with $\|z\|_1 \leq k$ can be written as a finite positive integer linear combination $\sum \alpha_i(g_i, g_i')$ with $\alpha_i \in \mathbb{Z}_{>0}$, $(g_i, g_i') \in \mathcal{G}_{\leq k} := \mathcal{G}_0 \cup \mathcal{G}_1 \cup \ldots \cup \mathcal{G}_k$, and $(g_i, g_i') \sqsubseteq_{j+1} (z, z')$ for all $i$.*

**Proof.** Let $(z, z') \in K_{j+1}^+ \cup K_{j+1}^-$ with $\|z\|_1 \leq k$. Clearly, by definition of $H_{j+1}^+ \cup H_{j+1}^-$, $(z, z')$ can be written as a finite positive integer linear combination $\sum \alpha_i(g_i, g_i')$ with $\alpha_i \in \mathbb{Z}_{>0}$, $(g_i, g_i') \in H_{j+1}^+ \cup H_{j+1}^-$, and $(g_i, g_i') \sqsubseteq_{j+1} (z, z')$ for all $i$. But $(g_i, g_i') \sqsubseteq_{j+1} (z, z')$ implies $\|g_i\|_1 \leq \|z\|_1 \leq k$, which, together with $(g_i, g_i') \in H_{j+1}^+ \cup H_{j+1}^-$, gives $(g_i, g_i') \in \mathcal{G}_{\leq k}$ and the proof is complete. $\square$

Note that it is crucial for the proof of the following corollary that, in case of $j < s$, we always choose the smallest non-negative integer $h'$ such that $(h, h') \in K_{j+1}$ when we form the input set $F$ to Algorithm 2.3. Moreover, note that $\mathcal{G}_0 = \emptyset$ if $j \geq s$ and that $\mathcal{G}_0 = \{\pi_{j+1}^n(p_{j+1}), -\pi_{j+1}^n(p_{j+1})\}$ if $j < s$.

**Corollary 3.2** *Let $(z, z') \in K_{j+1}^+ \cup K_{j+1}^-$, $\|z\|_1 = k+1$, and $|z'| < p_{j+1,j+1}$. Then there are precisely two possible situations:*

1. *There exists a vector $(v, v') \in \mathcal{G}_{\leq k}$ with $(v, v') \sqsubseteq (z, z')$. Then we already know in advance that $\text{normalForm}((z, z'), \mathcal{G}_{\leq k}) = 0$.*

2. *There does not exist a vector $(v, v') \in \mathcal{G}_{\leq k}$ with $(v, v') \sqsubseteq (z, z')$. Then $(z, z') \in \mathcal{G}_{k+1}$.*

**Proof.** First, suppose that there exists $(v, v') \in \mathcal{G}_{\leq k}$ with $(v, v') \sqsubseteq (z, z')$. Clearly, if $j < s$, $(v, v') \sqsubseteq_{j+1} (z, z')$ and $|z'| < p_{j+1,j+1}$ imply $(v, v') \notin \mathcal{G}_0 = \{\pi_{j+1}^n(p_{j+1}), -\pi_{j+1}^n(p_{j+1})\}$. Thus we have $v \neq 0$ and consequently $\|z - v\|_1 \leq k$. We conclude by Lemma 3.1 that there are finitely many (not necessarily different) elements $(g_i, g_i') \in \mathcal{G}_{\leq k}$ such that $(z - v, z' - v') = \sum (g_i, g_i')$ and $(g_i, g_i') \sqsubseteq_{j+1} (z - v, z' - v') \sqsubseteq_{j+1} (z, z')$ for all $i$. Choosing $(v, v')$ together with these $(g_i, g_i')$ in the algorithm normalForm we obtain $\text{normalForm}((z, z'), \mathcal{G}_{\leq k}) = 0$.

On the other hand, if there does not exist a vector $(v, v') \in \mathcal{G}_{\leq k}$ with $(v, v') \sqsubseteq_{j+1} (z, z')$, then $(z, z') \in \mathcal{G}_{k+1}$, that is, the vector $(z, z')$ cannot be written as a sum $(v_1, v_1') + (v_2, v_2')$ with $(v_i, v_i') \in (K_{j+1}^+ \cup K_{j+1}^-) \setminus \{0\}$ and $(v_i, v_i') \sqsubseteq_{j+1} (z, z')$, $i = 1, 2$. To see this, assume on the contrary that such vectors $(v_1, v_1')$ and $(v_2, v_2')$ do exist. Note that both vectors lie in the same orthant of $\mathbb{R}^{j+1}$ as $(z, z')$. From $(v_i, v_i') \sqsubseteq_{j+1} (z, z')$, $i = 1, 2$, and $|z'| < p_{j+1,j+1}$ we conclude that $v_i \neq 0$, $i = 1, 2$. Therefore, $\|v_i\|_1 \leq k$, $i = 1, 2$ which implies by Lemma 3.1 that $(v_1, v_1')$ and $(v_2, v_2')$ can both be written as finite positive integer linear combinations of elements from $\mathcal{G}_{\leq k}$ that all lie in the same orthant as $(v_1, v_1')$ and $(v_2, v_2')$. Put together, both combinations allow a finite positive integer linear representation of $(z, z')$ by elements from $\mathcal{G}_{\leq k}$ that all lie in the same orthant as $(z, z')$. But all the summands $(v, v') \in \mathcal{G}_{\leq k}$ in this representation of $(z, z')$ fulfill $(v, v') \sqsubseteq_{j+1} (z, z')$ in contradiction to our initial assumption that no such non-zero vector $(v, v')$ exists. Thus, $(z, z') \in \mathcal{G}_{k+1}$. □

In view of Corollary 3.2, considering the elements $(v, v')$ in $C$ by increasing norm $\|v\|_1$ leads to an algorithm to compute $H_{j+1}^+ \cup H_{j+1}^-$ that has the following advantages:

- The set $C$ of S-vectors need not be stored, as for $k = 2, 3, \ldots$, the elements in $\mathcal{G}_{k+1}$ are generated and immediately checked for reducibility within a simple loop. No superfluous vectors are stored.

- Thus, precisely $H_{j+1}^+ \cup H_{j+1}^-$ is computed.

- There is no reduction of vectors, only much cheaper reducibility tests are done.

- Upper bounds on variables can be easily be used directly to truncate the output set and to speed-up the computation.

It remains to give a criterion of when $H_{j+1}^+ \cup H_{j+1}^- = \mathcal{G}_{\leq k}$ that allows us to stop the inductive construction of the $\mathcal{G}_k$.

**Lemma 3.3** *Let $k \in \mathbb{Z}_{>0}$ satisfy $\mathcal{G}_{k+1} = \ldots = \mathcal{G}_{2k} = \emptyset$. Then $H_{j+1}^+ \cup H_{j+1}^- = \mathcal{G}_{\leq k}$.*

**Proof.** Note that the input set $F$ defined in Section 2 satisfies $F \subseteq \mathcal{G}_{\leq k}$. Therefore, by Lemma 2.2, also $\mathcal{G}_{\leq k}$ generates $K_{j+1}^+ \cup K_{j+1}^-$ over $\mathbb{Z}_+$. Thus, by Lemma 2.5, Algorithm 2.3 would return a set contain $H_{j+1}^+ \cup H_{j+1}^-$ if it were started with $\mathcal{G}_{\leq k}$ as input set. However, by Lemma 3.1 and Corollary 3.2, all S-vectors reduce to 0 with respect to $\mathcal{G}_{\leq k}$ and Algorithm 2.3 returns exactly $\mathcal{G}_{\leq k}$. As all elements in $\mathcal{G}_{\leq k}$ are minimal, we must have $\mathcal{G}_{\leq k} = H_{j+1}^+ \cup H_{j+1}^-$. □

For $j \geq s$, one question remains open from a practical point of view: As the remaining variables (with index bigger than $j$) can be freely permuted and then re-permuted at the end of the computation, what is a good selection strategy to choose the $(j+1)^{\text{st}}$ variable to be appended next? First computational experiments show that huge differences in sizes of intermediate sets $H_j^+$ and in running times do indeed occur.

## 4 Computation of Extremal Rays

Given a cone $\mathcal{C} = L \cap \mathbb{R}_+^n$, where $L$ is the linear space generated by $\{p_1, \ldots, p_s\} \subseteq \mathbb{Z}^n$. We may again assume that these generators have the same structure as in Section 2. Analogously to the computation presented in Section 2, let $\bar{K}_j := \{\pi_j^n(v) : v \in L\}$, $\bar{K}_j^+ := \bar{K}_j \cap (\mathbb{R}_+^{j-1} \times \mathbb{R}_+)$, $\bar{K}_j^- := \bar{K}_j \cap (\mathbb{R}_+^{j-1} \times \mathbb{R}_-)$. Note that $\bar{K}_j^+$ and $\bar{K}_j^-$ are pointed rational cones. Let $R_j^+$ and $R_j^-$ denote their sets of extreme rays. Again we will start with $R_1^+ = \{(p_{1,1})\}$, and compute $R_2^+$, $R_3^+$, $\ldots$, $R_n^+$ inductively. For each step $R_j^+ \to R_{j+1}^+$ we will again employ Algorithm 2.3. However, we have to use different specifications for the input set, the set of S-vectors, and the algorithm normalForm. Note that a similar speed-up as presented in Section 3 is possible for this algorithm, by choosing elements $(v, v')$ from $C$ by increasing value of $|\text{supp}(v)|$. We will, however, skip the details here.

The following Lemma can be proved in a similar way as Lemma 2.2.

**Lemma 4.1** *The set*

$$F = \begin{cases} \bigcup_{h \in R_j^+} \{(h, h') : (h, h') \in \bar{K}_{j+1}\} & \text{if } j \geq s, \text{ and} \\ \bigcup_{h \in R_j^+} \{(h, h') : (h, h') \in \bar{K}_{j+1}\} \cup \{\pi_{j+1}^n(p_{j+1}), -\pi_{j+1}^n(p_{j+1})\} & \text{if } j < s, \end{cases}$$

*generates* $\bar{K}_{j+1}^+ \cup \bar{K}_{j+1}^-$ *over* $\mathbb{R}_+$, *that is, every* $z \in \bar{K}_{j+1}^+ \cup \bar{K}_{j+1}^-$ *can be written as a non-negative linear combination of elements of* $F$.

Let $F$ as defined in Lemma 4.1 be the input set to the completion algorithm and claim the output to be a set $G$ which contains a positive scalar multiple for each element in $R_{j+1}^+ \cup R_{j+1}^-$. It remains to define the set S-vectors$(f, g)$ and the function normalForm$(s, G)$.

$$\text{S-vectors}((v, v'), (w, w')) := \begin{cases} \{(v - (v'/w')w, v' - (v'/w')w')\} & \text{if } v'w' < 0, \\ \emptyset & \text{otherwise.} \end{cases}$$

Note that all S-vectors have 0 as their last component. Behind the function normalForm$(s, G)$ there is the following algorithm.

**Algorithm 4.2** *(Normal Form Algorithm for Extreme Ray Computation)*
<u>Input:</u> *a vector $s$, a set $G$ of vectors*
<u>Output:</u> *a normal form of $s$ with respect to $G$*

<u>while</u> *there is some $g \in G$ such that* $\text{supp}(g) \subseteq \text{supp}(s)$ <u>do</u>
$\quad \alpha := \min\{s^{(i)}/g^{(i)} : i = 1, \ldots, j, g^{(i)} \neq 0\}$
$\quad s := s - \alpha g$
<u>return</u> *s*

Note that the first $j$ components of $s$ always remain non-negative during the normal form algorithm. Moreover, $\text{supp}(s - \alpha g) \subsetneq \text{supp}(s)$ and therefore, the normal form algorithm must always terminate.

**Lemma 4.3** *Algorithm 2.3 with the new specifications for the input set,* S-vectors *and* normalForm *(as presented in this section) terminates and returns a set containing* $R_{j+1}^+ \cup R_{j+1}^-$.

**Proof.** Termination of the algorithm is clear, since only a finite set $C$ of vectors is checked. Thus, it remains to show correctness of the algorithm.
By $G$ denote the set that is returned by the completion procedure and let $(z, z') \in R_{j+1}^+ \cup R_{j+1}^-$. By Lemma 4.2 and $F \subseteq G$, we can write $(z, z')$ as a positive linear combination $(z, z') = \sum \alpha_i(v_i, v'_i)$ for some $\alpha_i \in \mathbb{R}_{>0}$ and vectors $(v_i, v'_i) \in G$ with $\text{supp}(v_i) \subseteq \text{supp}(z)$ for all $i$. From the set of all such positive linear combinations $\sum \alpha_i(v_i, v'_i)$ choose one such that $\sum \alpha_i |v'_i|$ is minimal. In order to find such linear combination we have to solve the linear program

$$\min_{\alpha > 0} \{\sum \alpha_i |v'_i| : \sum \alpha_i(v_i, v'_i) = (z, z')\}$$

for given $(v_i, v'_i) \in G$ and $(z, z')$. This linear program is bounded from below by 0. Thus the minimal value is indeed attained by some choice of the $\alpha_i$.
Note that $\sum \alpha_i |v'_i| \geq |z'|$ with equality if and only if $z'$ and all $v'_i$ are either non-negative or non-positive.

If $\sum \alpha_i |v_i'| = |z'|$, then all $(v_i, v_i')$ lie in the same orthant as $(z, z')$ and we get $\mathrm{supp}((v_i, v_i')) \subseteq \mathrm{supp}((z, z'))$ for all $i$. Together with $(z, z') = \sum \alpha_i (v_i, v_i')$, $\alpha_i \in \mathbb{R}_{>0}$, and $(z, z') \in R_{j+1}^+ \cup R_{j+1}^-$ this representation must be trivial, that is, $(z, z') = \alpha_1 \cdot (v_1, v_1') \in G$, $\alpha_1 > 0$, and consequently, $G$ contains some positive scalar multiple of $(z, z')$ as claimed. Hence, we will assume on the contrary that $\sum \alpha_i |v_i'| > |z'|$ holds. Thus, there must exist $(v_{i_1}, v_{i_1}'), (v_{i_2}, v_{i_2}')$ such that $v_{i_1}' v_{i_2}' < 0$.

The sum $(v_{i_1}, v_{i_1}') - (v_{i_1}'/v_{i_2}')(v_{i_2}, v_{i_2}') = (v_{i_1} - (v_{i_1}'/v_{i_2}')v_{i_2}, 0)$ was reduced to 0 during the run of Algorithm 2.3 which gives a linear combination $(v_{i_1} - (v_{i_1}'/v_{i_2}')v_{i_2}, 0) = \sum \beta_k (w_k, w_k')$ for some positive scalars $\beta_k$ and some $(w_k, w_k') \in G$. Moreover, $\mathrm{supp}((w_k, w_k')) \subseteq \mathrm{supp}((v_{i_1} - (v_{i_1}'/v_{i_2}')v_{i_2}, 0))$ for all $k$, implying that

$$0 = \sum \beta_k |w_k'| = |\sum \beta_k w_k'| = |v_{i_1}' - (v_{i_1}'/v_{i_2}')v_{i_2}'| < |v_{i_1}'| - (v_{i_1}'/v_{i_2}')|v_{i_2}'|$$

holds. But if we now choose $\alpha = \min\{\alpha_{i_1}, \alpha_{i_2}/(-v_{i_1}'/v_{i_2}')\} = \min\{\alpha_{i_1}, -\alpha_{i_2} v_{i_2}'/v_{i_1}'\} > 0$, the representation

$$(z, z') = \sum \alpha \beta_k (w_k, w_k') + (\alpha_{i_1} - \alpha)(v_{i_1}, v_{i_1}') + (\alpha_{i_2} + \alpha(v_{i_1}'/v_{i_2}'))(v_{i_2}, v_{i_2}') + \sum_{i \neq i_1, i_2} \alpha_i (v_i, v_i')$$

contradicts the minimality of $\sum \alpha_i |v_i'|$. We conclude $\sum \alpha_i |v_i'| = |z'|$ and the claim follows. □

## 5 Applications

### 5.1 Decomposition of Vectors

**Problem.** Given $A \in \mathbb{Z}^{d \times n}$ and $u \in \mathbb{Z}_+^n$, find a representation $u = \sum_{i=1}^t \lambda_i h_i$ with $h_i \in H(\ker_\mathbb{Z}(A))$, $i = 1, \ldots, t$, and positive integral multipliers $\lambda_1, \ldots, \lambda_t$.

**Solution.** Compute a truncated Hilbert basis $\bar{H}(\ker_\mathbb{Z}(A))$ with upper bounds $u$ on the variables. Once $\bar{H}(\ker_\mathbb{Z}(A))$ has been computed, $u$ can be easily decomposed as a non-negative integer linear combination of elements of $\bar{H}(\ker_\mathbb{Z}(A))$ as follows: Find $v \in \bar{H}(\ker_\mathbb{Z}(A))$ such that $v \leq u$ and replace $u$ by $u - v$. Note that $u - v$ can also be written as a non-negative integer linear combination of elements of $\bar{H}(\ker_\mathbb{Z}(A))$ as $u - v \leq u$.

This solution approach can be sped-up once an element of some $H_j^+$ is found that extends to a vector $v \in \ker_\mathbb{Z}(A) \cap \mathbb{R}_+^n$ with $v \leq u$. Note that this vector $v$ must belong to $H_n^+$ as already its projection onto the first $j$ components is not decomposable. In case such a vector $v$ has been found, we may replace $u$ by $u - v$ and use the stronger upper bounds $u - v$ instead of $u$. Note that the computation of the truncated Hilbert basis need not be started again. It suffices to throw away unnecessary vectors from $H_j^+$ (that do not lie below the new upper bounds $u - v$) and to continue with the computation of $H_{j+1}^+$. If $H_n^+$ is reached, we can completely decompose the current $u$ or prove that it is not further decomposable and therefore $u \in \bar{H}(\ker_\mathbb{Z}(A)) \subseteq H(\ker_\mathbb{Z}(A))$.

### 5.2 Counting Magic Arrays and $t$-designs

**Problem.** Given a $d$ dimensional $n \times n \times \cdots \times n$ array with $n^d$ non-negative integer entries, wherein the sum of the entries in any coordinate (and in the main diagonal) direction equals a magic constant $s$. Find a formula for the number of magic arrays with given magic constant $s$.

**Start of Solution.** The method presented in [1] computes in a first step all non-decomposable magic arrays. This step is equivalent to the computation of the Hilbert basis of the cone of all magic arrays with non-negative real entries. It turned out that the computation of Hilbert bases is currently not the computational bottleneck in the method described in [1].
The computation of Hilbert bases is also employed in counting $t$-designs [6].

### 5.3 Integer Programming

**Problem.** Given $A \in \mathbb{Z}^{d \times n}$, $b \in \mathbb{Z}^d$, $c \in \mathbb{R}^n$, consider the integer program

$$\min_z \{c^\intercal z : Az = b, z \in \{0,1\}^n\}.$$

Given a 0-1 solution $z_0$ to $Az = b$, find a cheaper solution or prove optimality of $z_0$.
**Solution.** Transform the problem via $z^{(i)} \to 1 - (z')^{(i)}$ whenever $z_0^{(i)} = 1$, $z^{(i)} \to (z')^{(i)}$ otherwise, into a 0-1 problem with problem matrix $A'$, cost vector $c'$, and current feasible solution $z_0' = 0$. The above problem reduces to the computation of an element $v \leq 1$ in the truncated Hilbert basis of $\ker_\mathbb{Z}(A') \cap \mathbb{R}_+^n$ (with upper bounds 1) that has a strictly negative objective value $(c')^\intercal v$.

### 5.4 Dual cones

**Problem.** Given a cone $\mathcal{C}$ generated by integral vectors $p_1, \ldots, p_s$ in $\mathbb{R}^n$. Compute the generators $R(\mathcal{C}^D)$ and the minimal Hilbert basis $H(\mathcal{C}^D)$ of the dual cone

$$\mathcal{C}^D := \{v \in \mathbb{R}^n : z^\intercal v \geq 0, \forall z \in \mathcal{C}\} = \{v \in \mathbb{R}^n : p_i^\intercal v \geq 0, i = 1, \ldots, s\}.$$

Note the $R(\mathcal{C}^D)$ gives exactly the set of normal vectors of the facets of $\mathcal{C}$.
**Solution.** Let $I_s$ denote the $s \times s$ identity matrix and let $P$ be the $s \times n$ matrix that has $p_1^\intercal, \ldots, p_s^\intercal$ as rows. Moreover, let $L := \{(v,u) \in \mathbb{R}^{n+s} : Pv - I_s u = 0\}$ be generated over $\mathbb{R}$ by the vectors $(e_1, Pe_1), \ldots, (e_n, Pe_n)$. Note that these vectors also generate $\ker_\mathbb{Z}(P|-I_s)$ over $\mathbb{Z}$. Therefore, we have $\mathcal{C}^D = \pi_n^{n+s}(L \cap (\mathbb{R}^n \times \mathbb{R}_+^s))$ and thus

$$R(\mathcal{C}^D) = \{v \in \mathbb{R}^n : u = Pv \in R(P^\intercal)\} \text{ and } H(\mathcal{C}^D) = \{v \in \mathbb{Z}^n : u = Pv \in H(P^\intercal)\},$$

where $R(P^\intercal)$ and $H(P^\intercal)$ denote the set of extreme rays and the Hilbert basis of the cone which is the intersection of $\mathbb{R}_+^s$ with the subspace of $\mathbb{R}^s$ spanned over $\mathbb{R}$ by the columns of $P$ (rows of $P^\intercal$).

### 5.5 Hilbert Bases of Cones given by their Generators

**Problem.** Given a cone $\mathcal{C}$ generated by integral vectors $p_1, \ldots, p_s$ in $\mathbb{R}^n$, compute the minimal Hilbert basis $H(\mathcal{C})$ of $\mathcal{C}$.
**Solution.** Since $\mathcal{C} = (\mathcal{C}^D)^D$, this problem can be solved by first computing the generators $R(\mathcal{C}^D)$ of the dual cone and by then computing the Hilbert basis of the dual of the dual cone $\mathcal{C}^D$, which is the Hilbert basis of $\mathcal{C}$.

## 6 Computational Experience

The algorithm presented in Sections 2 and 3 were implemented into a computer program that we called MLP. It is written in C and downloadable from http://www.testsets.de. Besides (truncated) Graver test sets [8], MLP allows the computation of a (truncated) set of minimal generators of the monoid $\Lambda \cap \mathbb{R}_+^n$ for any lattice $\Lambda \subseteq \mathbb{Z}^n$, in particular for $\Lambda = \ker_\mathbb{Z}(A)$.

In [1], MLP was used to solve some non-trivial instances. For example, it took about 10 days on a 1GHz PC with 4GB RAM running linux to compute the set of all minimal $6 \times 6$ magic squares, yielding a Hilbert basis containing $522,347$ elements. These examples can be found on the benchmark section of http://www.testsets.de.

**Acknowledgments.** The author wants to thank Jesús De Loera, UC Davis, for many helpful discussions and remarks. Moreover, the author gratefully acknowledges the support by AFOSR Grant F49620-01-0327 and by NSF Grant DMS-0073815.